\documentclass[12pt,a4paper]{article}

\usepackage{amsmath,amstext,amssymb,amscd,euscript}

\usepackage[english,russian]{babel}
\usepackage[cp1251]{inputenc}

\newcommand{\Xcomment}[1]{}

\newtheorem{theorem}{Теорема}[section]
\newtheorem{lemma}[theorem]{Лемма}
\newtheorem{corollary}[theorem]{Следствие}

\newtheorem{prop}[theorem]{Предложение}

\makeatletter \@addtoreset{equation}{section} \makeatother

\newenvironment{proof}{\noindent{\bf Доказательство}\/}%
{\hfill$\qed$\medskip}

\def\qed{\ \vrule width.1cm height.3cm depth0cm}

\newenvironment{numitem1}{\refstepcounter{equation}\begin{enumerate}%
\item[(\thesection.\arabic{equation})]}{\end{enumerate}}

\newcommand{\refeq}[1]{(\ref{eq:#1})}  

 \makeatletter
\renewcommand{\section}{\@startsection{section}{1}{0pt}%
{-3.5ex plus -1ex minus -.2ex}{2.3ex plus .2ex}%
{\normalfont\Large}}
 \makeatother

 \makeatletter
\renewcommand{\subsection}{\@startsection{subsection}{2}{0pt}%
{-3.0ex plus -1ex minus -.2ex}{1.5ex plus .2ex}%
{\normalfont\normalsize\bf}}
 \makeatother

 \newcommand{\SEC}[1]{\ref{sec:#1}}  

\def\Rset{{\mathbb R}}

\def\Zset{{\mathbb Z}}

\def\Cscr{{\cal C}}

\def\Tscr{{\cal T}}

\def\tilde{\widetilde}
\def\hat{\widehat}

\def\eps{\varepsilon}

\def\deltain{\delta^{\rm in}}
\def\deltaout{\delta^{\rm out}}

\def\excess{{\rm ex}}
\def\Excess{{\rm Old}}

\def\New{{\rm New}}
\def\valu{{\rm val}}

\begin{document}

 \title{On stable flows and preflows}

 \author{Alexander V.~Karzanov
\thanks{Central Institute of Economics and Mathematics of
the RAS, 47, Nakhimovskii Prospect, 117418 Moscow, Russia; email:
akarzanov7@gmail.com.}
 }

\date{}

 \maketitle

\begin{quote}
\textbf{Abstract.} In 2010s Fleiner introduced a notion of stable flows in
directed networks and showed that such a flow always exists and can be found by
use of a reduction to the stable allocation problem due to Baiou and Balinski.
Recently Cseh and Matuschke devised a direct strongly polynomial algorithm. In
this paper we give an alternative algorithm to find a stable flow in a network
with several sources and sinks. It is based on an idea of \emph{preflows}
(appeared  in 1970s in a faster algorithm for the classical max-flow problem),
and runs in $O(nm)$ time for a network with $n$ vertices and $m$ edges. The
results are further generalized to a larger class of objects, so-called stable
\emph{quasi-flows} with bounded excesses in non-terminal vertices. (The paper
is written in Russian.)
\smallskip

\noindent\emph{Keyworks}:\,  Stable flow in a network, Stable allocation, Preflow, Quasi-flow
\end{quote}


\baselineskip=15pt
\parskip=2pt

\section{Введение}  \label{sec:intr}

Область теоретических и прикладных задач о стабильных договорах явилась предметом интенсивного изучения в математической экономике, теории игр и комбинаторной оптимизации, и ей посвящены многочисленные работы нескольких последних десятилетий. Отправной точкой многих исследований послужила классическая работа Гейла и Шепли~\cite{GS} о \emph{стабильных марьяжах} (сокращенно SM).

Согласно одной из распространенных формулировок этой задачи имеется двудольный граф $G=(V,E)$, и для каждой вершины $v$ задан (строгий) линейный порядок $<_v$ на инцидентных ребрах. (В~\cite{GS} рассматриваются полные двудольные графы, но это не существенно. В популярной интерпретации ребра в $G$ представляют возможные союзы персон разного пола, а порядок $<_v$ -- предпочтения персоны $v$: если для ребер $e=vu$ и $e'=vw$ выполняется $e<_v e'$, то $v$ предпочитает союз с $u$ союзу с $w$.) В задаче требуется найти паросочетание (matching) $M\subseteq E$, являющееся стабильным относительно всех этих порядков. Это означает, что для любого ребра $e$ из $E-M$ есть ребро $e'\in M$ такое, что $e$ и $e'$ имеют общую вершину $v$, и при этом выполняется $e'<_v e$. Было показано, что стабильное паросочетание в двудольном графе всегда существует, и оно может быть построено комбинаторным алгоритмом с линейной верхней оценкой числа действий (\emph{временн\'ой сложностью}) $O(n+m)$, где n и m -- число вершин и ребер в $G$, соответственно.

В последующих работах многих авторов были исследованы различные обобщения задачи SM. Укажем два направления обобщений, связанные с графами (оставляя в стороне постановки, в которых в договорах могут участвовать более двух агентов, или предпочтения агентов задаются другими способами). Одно из них -- переход от двудольного к произвольному графу $G$. Соответствующий аналог задачи SM, известный под названием задачи о \emph{стабильном размещении пар соседей по комнатам} (stable roommates problem), был исследован Ирвингом в работе~\cite{irv}, где был дан алгоритм линейной сложности для построения стабильного паросочетания в $G$ либо доказательства его несуществования. Важные дополнительные структурные и алгоритмические результаты были представлены в~\cite{tan}.

Другой тип обобщений, который более важен для нас, оставляет граф $G=(V,E)$ двудольным, но  добавляет числовые параметры. Среди задач этого типа весьма общей выглядит задача о \emph{стабильном распределении} (stable allocation problem), сокращенно SA, введенная Бейу и Балински~\cite{BB}. Здесь распределением считается назначение каждому ребру $e\in E$ величины $x(e)\ge 0$, не превышающей предписанной пропускной способности $c(e)$, и при этом сумма назначений по ребрам, инцидентным вершине $v\in V$, не должна превышать предписанной ``квоты'' $q(v)$. (В случае задачи SM, все $c(e)$ и $q(v)$ равны единице, а $x(e)$ принимает значение 0 или 1. В общем случае задачи SA, число $x(e)$ на ребре $e=uv$ может пониматься, например, как размер участия ``работника'' $u$ в ``работе'' $v$.) В~\cite{BB} доказывается разрешимость SA при любых неотрицательных вещественных $c,q$ (давая целочисленное $x$ при целочисленных $c,q$)  и предложен сильно полиномиальный алгоритм решения, т.е. имеющий временн\'ую сложность, зависящую только от размеров графа и выражающуюся полиномом от $n$ и $m$.
Дин и Мунши~\cite{DM} описали улучшенную версию алгоритма из~\cite{BB}, которая строит решение за время $O(nm)$; более того, они показали, что можно добиться даже временн\'ой оценки $O(m\,\log n)$, если применить для ряда процедур мощные структуры данных, такие как динамические и самонастраивающиеся деревья Слейтора и Таржана~\cite{ST1,ST2}.
(Это дает теоретическое ускорение, но алгоритм такого рода громоздкий и едва ли может быть применен для практических целей.)

В свою очередь задача SA может быть представлена как частный случай задачи о \emph{стабильном потоке} (сокращенно SF). Последняя была сформулирована Флейнером в~\cite{flein} (как обобщение задачи Островского~\cite{ostr} для ациклических сетей с единичными пропускными способностями ребер). В ней задана сеть, состоящая из ориентированного графа $G=(V,E)$ с пропускными способностями $c(e)\ge 0$ ребер $e\in E$ и двумя выделенными вершинами (``терминалами'') $s$ и $t$. Для каждой \emph{внутренней} вершины $v\in V-\{s,t\}$ заданы линейный порядок на входящих ребрах и линейный порядок на выходящих ребрах. (Допустимый) \emph{поток} -- это неотрицательная вещественная функция  $f$ на ребрах, удовлетворяющая верхним ограничениям по пропускным способностям и имеющая нулевые эксцессы для всех внутренних вершин, где под эксцессом в вершине $v$ понимается разность $\excess_f(v)$ между общим потоком по входящим ребрам и общим потоков по выходящим ребрам в $v$. (Внутреннюю вершину можно интерпретировать как ``игрока'', ``торговца'' или ``агента'', который, получая некоторый объем  однородного продукта по входящим ребрам, посылает его далее по выходящим ребрам, руководствуясь своей функцией полезности, зависящей от указанных порядков на ребрах.) Поток считается стабильным, если он не допускает локальных улучшений, использующих ``ненасыщенные'' пути; точное определение будет дано а разделе~\SEC{defin}.

Между задачами SA и SF есть тесная связь. Задача SA с двудольным графом $G=(V_1\sqcup V_2, \,E)$ сводится к SF с графом, получаемым из $G$ (с ребрами, ориентированными от $V_1$ к $V_2$) путем добавления терминалов $s$ и $t$, ребер $(s,u)$ с пропускными способностями $q(u)$ для вершин $u$ из доли $V_1$, и ребер $(v,t)$ с пропускными способностями $q(v)$ для вершин $v$ из доли $V_2$.

С другой стороны, Флейнер~\cite{flein} показал, что задача SF с графом $G=(V,E)$ может быть сведена к задаче SA с графом, получаемым расщеплением вершин в $G$ и добавлением $O(|V|)$ новых ребер. В результате было установлено существование стабильного потока для любой сети (и целочисленного стабильного потока при целочисленном $c$), и возможность построения такого потока при помощи алгоритмов для SA, сохраняя аналогичную временн\'ую сложность.

Впоследствии появились прямые алгоритмы нахождения стабильного потока. Недавно Чех и Матушке~\cite{CM} предложили прямой алгоритм для сети с одним источником и одним стоком, который имеет сложность $O(nm)$ и основан на комбинации идей метода увеличивающих путей Форда и Фалкерсона для задачи о максимальном потоке и метода ``отсроченных принятий'' (deferred aceptance), восходящего к Гейлу и Шепли~\cite{GS}.

В настоящей работе предлагается альтернативный алгоритм построения стабильного потока в сети $N=(G,S,T,c)$ с произвольными множествами источников $S$ и стоков $T$ при условии отсутствия ребер входящих в источники или выходящих из стоков. Алгоритм прямой и чисто комбинаторный (не апеллирует к задаче SA и не использует сложных структур данных), он основан на методе предпотоков. Напомним, что \emph{предпотоком} в сети называется неотрицательная функция на ребрах, ограниченная пропускными способностями и имеющая неотрицательные эксцессы во всех внутренних вершинах. (Это понятие было введено в~\cite{karz} и использовалось в алгоритме нахождения максимального потока в сети.) Заметим, что предпотоки уже применялись ранее в работе~\cite{CMS} для построения стабильного потока в сети с целочисленными пропускными способностями за псевдополиномиальное время (линейно зависящее от суммы пропускных способностей ребер).

Алгоритм имеет базовую и модифицированную (ускоренную) версии. Обе начинаются с построения некоторого стабильного предпотока, и на каждой последующей итерации текущий стабильный предпоток перестраивается с целью избавления от положительных эксцессов во внутренних вершинах. Как только эксцессы всех внутренних вершин обнулятся, будет построен искомый стабильный поток (причем целочисленный при целочисленном $c$). Базовый алгоритм конечен при любых неотрицательных вещественных пропускных способностях $c$. Модифированный алгоритм сильно полиномиальный; он использует дополнительные преобразования предпотоков и строит стабильный поток с временн\'ой сложностью $O(nm)$ (подобно алгоритму в~\cite{CM}).

Мы затем рассматриваем более общую задачу для сети $N=(G,S,T,c)$, в которой для каждой внутренней вершины $v$ заданы два параметра $\beta(v)\ge 0$ и $\gamma(v)\ge 0$ и требуется найти стабильный ``квазипоток'' $f$ с ограничениями вида $-\beta(v)\le \excess_f(v)\le \gamma(v)$. (Это превращается в стабильный поток при $\beta,\gamma=0$.) Показывается, что такой ``квазипоток'' существует и также может быть найден за время $O(nm)$. (В приложениях число $\gamma(v)$ можно понимать как разрешение ``агенту'' $v$ распоряжаться по своему усмотрению частью полученного продукта в размере, не превышающем $\gamma(v)$, а число $\beta(v)$ -- привлекать ``со стороны'' дополнительный продукт в размере не более $\beta(v)$.)

Данная статья организована следующим образом. В разделе~\SEC{defin} приводятся основные определения и точные формулировки задач SF и SA. В разделе~\SEC{bas_alg} излагается базовый алгоритм нахождения стабильного потока в ациклической сети методом предпотоков и показывается его конечная сходимость (Предложение~\ref{pr:finite}). Модифицированная версия алгоритма, имеющая временн\'ую сложность $O(nm)$, описывается в разделе~\SEC{modif_alg}. Раздел~\SEC{general} содержит обобщения полученных результатов на предпотоки и квазипотоки с ограниченными эксцессами (Предложения~\ref{pr:gamma_pre} и~\ref{pr:beta-gamma}). В заключительном разделе~\SEC{concl} обсуждаются три дополнительных свойства: 1) максимум величин (т.е. суммарных эксцессов в стоках) для стабильных предпотоков достигается на стабильном потоке; 2) стабильные потоки в фиксированной сети имеют одинаковые величины; и 3) стабильные потоки образуют решетку. (Свойства в пунктах 2 и 3 показываются в~\cite{flein} через сведение к соответствующим свойствам в задаче SA; мы описываем прямые конструкции.)


\section{Определения и постановки} \label{sec:defin}

Мы рассматриваем сеть $N=(G,S,T,c)$, состоящую из \emph{ориентированного} графа $G=(V,E)$ (без петель и кратных ребер), выделенных непересекающихся подмножеств вершин $S$ (\emph{источники}) и $T$ (\emph{стоки}), называемых также \emph{терминалами},  и функции $c:E\to\Rset_+$  \emph{пропускных способностей} ребер. (Здесь и далее $\Rset_+$ и $\Zset_+$ -- множества неотрицательных вещественных и неотрицательных целых чисел, соответственно.) Для вершины $v\in V$ обозначим $\deltain(v)$ и $\deltaout(v)$ множества ребер \emph{входящих} в $v$ и \emph{выходящих} из $v$, соответственно. Для большей простоты нашего изложения мы предполагаем, что

\begin{numitem1} \label{eq:P}
~$\deltain(s)=\emptyset$ для всех $s\in S$, и $\deltaout(t)=\emptyset$ для всех $t\in T$.
 \end{numitem1}

\noindent\textbf{Определения.} Функция $f:E\to\Rset_+$ называется \emph{допустимой} (относительно $c$), если она удовлетворяет ограничениям $f(e)\le c(e)$ для всех ребер $e\in E$. Определим \emph{эксцесс} для $f$ в вершине $v\in V$ как
  \begin{equation} \label{eq:excess}
  \excess_f(v):=\sum\nolimits_{e\in\deltain(v)} f(e) -\sum\nolimits_{e\in\deltaout(v)} f(e).
  \end{equation}
Допустимая функция $f$ называется \emph{предпотоком} в $N$ (следуя терминологии в~\cite{karz}), если каждая вершина $v\in V-S$ имеет  неотрицательный эксцесс: $\excess_f(v)\ge 0$. \emph{Поток} (из $S$ в $T$) -- это предпоток $f$ с нулевыми эксцессами всех \emph{внутренних} (нетерминальных) вершин $v\in V-(S\cup T)$, а его \emph{величиной} $\valu(f)$ считается $\excess_f(T):=\sum_{t\in T}\,\excess_f(t)$.
  \medskip

\emph{Путь} в $G$ -- это последовательность $P=(v_0,e_1,v_1,\ldots,e_k,v_k)$,  где $e_i$ -- ребро, соединяющее вершины $v_{i-1}$ и $v_i$. Ребро $e_i$ в $P$ называется  \emph{прямым} (forward), если $e_i=v_{i-1}v_i$, и \emph{обратным} (backward), если $e_i=v_iv_{i-1}$. (Мы обозначаем ребро, выходящее из $u$ и входящее в $v$ как $uv$, вместо обычного $(u,v)$). Путь называется \emph{ориентированным} (directed), если все его ребра прямые, и называется \emph{простым} (simple), если все его вершины различны. Противоположный (обратный) путь $(v_k,e_k,v_{k-1},\ldots,e_1,v_0)$ обозначается $P^{-1}$. Путь из вершины $u$ в вершину $v$ может быть назван $u$--$v$ \emph{путем}. Если не оговорено противное, то, говоря о пути,  мы считаем его \emph{нетривиальным}, т.е. имеющим по крайней мере одно ребро.

Для допустимой функции $f$ ребро $e$ с $f(e)=c(e)$ ($f(e)<c(e)$; $f(e)=0$) называется \emph{насыщенным} (соответственно, \emph{ненасыщенным}; \emph{свободным} (от $f$)). Путь считается ненасыщенным, если таковы все его ребра.

В рассматриваемой нами задаче каждая \emph{внутренняя} вершина $v$ сети $N$ снабжена линейным (полным строгим) порядком $<^-_v$ на множестве $\deltain(v)$ и линейным порядком $<^+_v$ на множестве $\deltaout(v)$. Они интерпретируются как ``отношения предпочтения'', а именно,  $e<^-_v e'$ означает, что вершина (``агент'') $v$ предпочитает ребро $e$ ребру $e'$; в этом случае мы также будем говорить, что $e$ расположена в $\deltain(v)$ \emph{раньше} или \emph{левее} $e'$, а $e'$ -- \emph{позже} или \emph{правее} $e$. В соответствии с порядком $<^-_v$ множество $\deltain(v)$ организуется в виде списка (double-linked list); таким образом, первый (последний) элемент в списке наиболее (соответственно, наименее) предпочтительный. И аналогично для порядка $<^+_v$.
 \medskip

\noindent\textbf{Определение.} Поток $f$ в сети $N$ с указанными порядковыми оснащениями для внутренних вершин называется \emph{стабильным}, если каждый \emph{ненасыщенный ориентированный} путь $P=(v_0,e_1,v_1,\ldots,e_k,v_k)$ удовлетворяет по крайней мере одному из  следующих двух условий:
 \begin{numitem1} \label{eq:v0}
начальная вершина $v_0$ внутренняя, и $P$  \emph{доминируется} в $v_0$; это означает, что каждое ребро $e\in\deltaout(v_0)$ позднее $e_1$ свободно от $f$;
  \end{numitem1}
 \begin{numitem1} \label{eq:vk}
конечная вершина $v_k$ внутренняя, и $P$  \emph{доминируется} в $v_k$, т.е. каждое ребро $e\in\deltain(v_k)$ позднее $e_k$ свободно от $f$
  \end{numitem1}
(допуская случай $v_0=v_k$). В частности, нет ненасыщенного ориентированного пути из $S$ в $T$. Предпоток $f$ называется стабильным в аналогичном случае.

Нас интересует \emph{задача о стабильном потоке (SF)}: найти стабильный поток для сети $N$. Флейнер~\cite{flein} установил, что эта задача всегда разрешима; более того:

 \begin{numitem1} \label{eq:flein}
если сеть $N=(G=(V,E),S,T,c)$ целочисленная (т.е. $c\in \Zset_+^E$), то существует целочисленный стабильный поток.
  \end{numitem1}

Это доказывалось путем редукции задачи SF для $N$ к \emph{задаче о стабильном распределении (SA)} на двудольном графе $G'=(V',E')$ с весами $c'(e)\ge 0$ ребер $e\in E'$ и ``квотами'' $q(v)\ge 0$ вершин $v\in V'$. (Строго говоря, \cite{flein} рассматривал случай, когда $|S\cup T|=2$, и условие~\refeq{P} не накладывается, но указанные свойства верны и для нашего случая.)  При этом редукция линейна по размерам, а именно:  $|V'|=2|V|$ и $|E'|<|E|+2|V|$, и сохраняет целочисленность: $c'$ и $q$ целочисленные при целочисленном $c$. Задача SA была введена и исследована в~\cite{BB}, где были доказаны ее разрешимость для любой сети и свойство целочисленности и предложен сильно полиномиальный алгоритм решения.

Как было отмечено во Введении, в~\cite{DM} показано, что при использовании продвинутых структур данных, т.н. динамических и самонастраивающихся деревьев, задачу о стабильном распределении в графе с $n$ вершинами и $m$ ребрами можно решить с временн\'ой сложностью  $O(m \log n)$ (а без использования таких структур можно добиться оценки $O(nm)$). Это, в силу указанной редукции, приводит к аналогичной оценке алгоритмической сложности и для задачи SF. Прямой алгоритм для SF в случае $|S|=|T|=1$, предложенный в~\cite{CM}, имеет сложность $O(nm)$.

В настоящей работе предлагается альтернативный прямой алгоритм нахождения стабильного потока в произвольной сети  $N=(G,S,T,c)$ (при условии~\refeq{P}), основанный на методе предпотоков. В следующем разделе мы описываем базовую версию алгоритма (конечную при любом $c$), а в разделе~\SEC{modif_alg} -- модифицированную версию (с временн\'ой сложностью $O(nm)$).

  \section{Базовый алгоритм} \label{sec:bas_alg}

Алгоритм нахождения стабильного потока в сети $N=(G=(V,E),S,T,c)$ состоит из последовательности итераций; как правило (но не всегда) итерация имеет две \emph{фазы}:  \emph{балансирование} (balancing) и \emph{достройка} (pushing)  (подобно структуре этапа (большой итерации) в алгоритме построения максимального потока в~\cite{karz}). Каждая итерация преобразует один блокирующий предпоток в другой, и процесс завершается, когда текущий предпоток становится потоком. Подчеркнем, что термин ``блокирующий'' заимствован из языка в ~\cite{karz} и имеет иной смысл, чем тот, что обычно понимают в задачах о стабильности.
  \medskip

\noindent\textbf{Определения.} Для предпотока $f$ в сети $N$ скажем, что вершина $v$ \emph{эксцессная}, если $\excess_f(v)>0$. Предпоток $f$
называется \emph{блокирующим}, если каждый ориентированный путь из $S$ в $T$ имеет насыщенное ребро. Если к тому же всякий ориентированный путь, идущий из эксцессной внутренней вершины в $T$ имеет насыщенное ребро, то мы называем $f$ \emph{полностью блокирующим}.
  \medskip

\noindent\textbf{Замечание 1.} Полностью блокирующий предпоток не обязательно стабильный, и наоборот. В то же время всякий стабильный \emph{поток} является блокирующим (и автоматически полностью блокирующим). В целом уже в ранжированной ациклической сети построение стабильного потока выглядит более сложной задачей, чем построение блокирующего потока. Последняя решается на этапе алгоритма в~\cite{karz} за время $O(n^2)$, что быстрее, чем $O(nm)$ в алгоритме раздела~\SEC{modif_alg} для SF.
 \medskip

\noindent\textbf{Начальная итерация.} На этой и последующих итерациях поддерживаются два списка $\Excess$ и $\New$, элементами которых являются эксцессные внутренние вершины текущуго предпотока. Также для каждой внутренней вершины $v$ указан элемент $\tilde e(v)$, который либо принимает пустое значение: $\tilde e(v)=\{\emptyset\}$, либо является выделенным ненасыщенным ребром в $\deltaout(v)$, называемым \emph{активным} ребром. В начале полагается $\Excess:=\New:=\emptyset$, и для каждого $v\in V-(S\cup T)$ активным полагается первое (самое предпочтительное) ребро в $\deltaout(v)$.

Начальная итерация состоит из одной фазы -- достройки, которая начинается с тривиального блокирующего предпотока $f$, определяемого как $f(e):=c(e)$ для всех ребер $e\in\deltaout(s)$, $s\in S$, и $f(e):=0$ для остальных ребер $e$. Каждая вершина $v$, соединенная ребром $e=sv$ с источником $s\in S$, заносится в список $\New$ (так как эксцесс в $v$ становится положительным). Затем мы сканируем вершины текущего списка $\New$.

Для очередной вершины $v\in \New$  мы изменяем $f$ на $\deltaout(v)$ так, чтобы либо свести эксцесс в $v$ к нулю, либо насытить все ребра в $\deltaout(v)$ (либо и то и другое). Более точно, ребра $e\in\deltaout(v)$ обрабатываются в порядке $<^+_v$ (``слева-направо''), начиная с активного ребра $e=\tilde e(v)$. Если текущий эксцесс $\Delta:=\excess_f(v)$ еще ненулевой, назначаем $f(e):=\min\{c(e),\Delta\}$. Если новый эксцесс обнулился (в случае $c(e)\ge\Delta$), обработка $v$ заканчивается; иначе (когда $c(e)<\Delta$), переходим к следующему ребру в списке $\deltaout(v)$, и т.д. Также при изменении $f$ на ребре $e=vw$, если вершина $w$ (в которой происходит увеличение эксцесса) не фигурирует ни в $\Excess$, ни в $\New$, то мы добавляем $w$ в текущий список $\New$.

При окончании работы с $v$ эта вершина удаляется (быть может временно) из списка $\New$, и если эксцесс в $v$ все еще ненулевой, $v$ заносится в список $\Excess$. Новым активным ребром $\tilde e(v)$ становится самое левое ненасыщенное ребро, а если такого нет  (в частности, когда $\excess_f(v)>0$), то полагается $\tilde e(v):=\{\emptyset\}$. Затем переходим к обработке другой вершины в $\New$, и т.д.. Итерация заканчивается, когда $\New$ становится пустым.
 \smallskip

Заметим, что в процессе итерации одна и та же внутренняя вершина $v$ может несколько раз появляться и исчезать в $\New$. Тем не менее, мы увидим позднее, что начальная (и каждая последующая) итерация всегда заканчивается. Можно видеть, что:
  \begin{numitem1} \label{eq:initial}
При завершении начальной итерации полученная функция $f$ является полностью блокирующим предпотоком $f$ с тем свойством, что для каждой внутренней вершины $v$: если $\tilde e(v)=\{\emptyset\}$, то все ребра в $\deltaout(v)$ насыщены, а если $\tilde e(v)\ne\{\emptyset\}$, то  все ребра в $\deltaout(v)$ раньше (левее) $e$ насыщены, а все ребра после $e$ свободны от $f$. При этом все эксцессные внутренние вершины $v$ (и только они) включены в список $\Excess$, и для таких $v$ выполняется $\tilde e(v)=\{\emptyset\}$.
 \end{numitem1}

\noindent Из~\refeq{initial} легко заключить выполнение~\refeq{v0} для всех ненасыщенных ориентированных путей, и следовательно, начальный предпоток $f$ стабильный.

В оставшейся части этого раздела мы сначала описываем фазу балансирования для общей итерации. Затем уточняем условия, которым должен удовлетворять текущий предпоток перед началом фазы достройки. Наконец, мы описываем эту фазу для общей итерации.
 \medskip

\noindent\textbf{Балансирование.} ~Эта процедура производится, когда множество $\Excess$ эксцессных внутренних вершин для текущего предпотока $f$ непусто (в то время как $\New=\emptyset$).  Она применяется к одной выбранной вершине $v\in\Excess$. Положим $\Delta:=\excess_f(v)$ ($>0$), и пусть $e$ -- последнее (самое правое, наименее предпочтительное) ребро в $\deltain(v)$ с $f(e)>0$. Уменьшим $f$ в $e$ на $\min\{\Delta,f(e)\}$. Если эксцесс в $v$ стал нулевым (в случае $\Delta\le f(e)$), заканчиваем. Иначе (когда $\Delta>f(e)$), берем последнее ребро $e'$ с $f(e')>0$ для обновленного $f$ и уменьшаем $f$ в $e'$ аналогичным образом. И так далее, пока эксцесс в $v$ не станет нулевым.

Для каждого ребра $e=uv$, в котором при данном балансировании было уменьшение $f$, проверяем вершину $u$, и если она внутренняя и не содержится в $\Excess$, то добавляем ее в список $\New$. (Таким образом, $\New$ -- это множество новых эксцессных вершин $u$, образовавшихся в результате уменьшения $f$ на ребрах $uv$.)  Cамое левое (хронологически последнее) ребро, где было уменьшение $f$, назовем \emph{критическим} и обозначим $\hat e=\hat e(v)$. Это ребро и все ребра $e=uv$ в $\deltain(v)$ после (правее) него помечаются как  \emph{закрытые} (только для последующего увеличения!). При этом если такое $e$ было активным ребром в $\deltaout(u)$, то новым активным ребром назначается первое после $e$ незакрытое  (при более ранних балансированиях) ребро в $\deltaout(u)$, а если такого нет, то полагается $\tilde e(u):=\{\emptyset\}$.
 \medskip

\noindent\textbf{Замечание 2.} Множества $\Excess$ и $\New$ задаются в виде обычных двухсторонних списков. В процессе алгоритма надо поддерживать величины эксцессов для внутренних вершин (корректируя их за $O(1)$ действий при каждом изменении $f$).
\medskip

Мы предполагаем (по индукции), что после балансирования в вершине $v$ текущее $f$ является предпотоком, который не обязан быть полностью блокирующим, но для которого выполнены следующие свойства:

 \begin{itemize}
\item[(C1)] Каждое ненасыщенное ребро в $\deltaout(s)$, $s\in S$, закрытое.
\item[(C2)] Если для эксцессной внутренней вершины $u$ множество $\deltaout(u)$ содержит незакрытое ненасыщенное ребро, то $u\in\New$.
\item[(C3)] Для каждой внутренней вершины $u$: (a)  если  в $\deltaout(u)$ есть активное ребро $\tilde e(u)\ne\{\emptyset\}$, то это ребро ненасыщенное и незакрытое, все ребра после него свободны от $f$, а каждое ребро до $\tilde e(u)$ -- насыщенное или закрытое (ставшее таковым при последнем или более раннем балансировании); (b) если $\tilde e(u)=\{\emptyset\}$, то все ребра в $\deltaout(u)$ насыщенные или закрытые; и (c) если $u$ когда-либо подвергалась балансированию, то $\tilde e(u)=\{\emptyset\}$, в свою очередь в $\deltain(u)$ есть критическое ребро $\hat e(u)$, это ребро ненасыщенное, а все ребра после него свободны от $f$.
\end{itemize}

\noindent\textbf{Достройка} (после балансирования в вершине $v$). Она состоит в увеличении текущего предпотока $f$ на некоторых ребрах, чтобы сделать его полностью блокирующим; это схоже с построением начального предпотока, но имеет ряд особенностей. Достройка немедленно заканчивается, когда начальный список $\New$ пустой. Иначе она начинается с некоторой вершины $u\in\New$, и мы стремимся максимально уменьшить эксцесс в $u$. Для этого сканируем ребра в $\deltaout(u)$ в порядке $<_u^+$, пропуская закрытые ребра. Сканирование начинается с активного ребра $\tilde e(u)\ne \{\emptyset\}$ (если $\tilde e(u)= \{\emptyset\}$, то вершина $u$ просто переносится из $\New$ в $\Excess$, и работа с $u$ заканчивается). Как и при начальной итерации, для очередного обрабатываемого ребра $e=uw$ полагаем $f(e):=\min\{c(e),\Delta\}$, где $\Delta$ - текущий эксцесс в $u$, и одновременно добавляем вершину $w$ в текущее множество $\New$, если его еще не было в $\Excess\cup\New$ (так как эксцесс в $w$ увеличивается). Если новый эксцесс в $u$ все еще ненулевой, переходим к следующему незакрытому ребру в $\deltaout(u)$, и т.д. Если эксцесс обнулился, то работу с $u$ заканчиваем и удаляем $u$  из $\New$. Если все ребра отсканированы, но эксцесс в $u$ все еще ненулевой, то переносим $u$ из $\New$ в $\Excess$. Окончив работу с $u$ и соответственно скорректировав активное ребро $\tilde e(u)$ в $\deltaout(u)$ (получая $\tilde e(u)$ как в~(C3)(a),(b)), переходим к другой вершине в текущем $\New$, и т.д. Заканчиваем процесс достройки, когда множество $\New$ становится пустым.
\medskip

При завершении достройки (которая конечна согласно Предложению~\ref{pr:finite} ниже) построенное $f$ является предпотоком, и можно видеть, что $f$ имеет свойства (C1) и (C3), а (C2) заменяется на следующее (ср.~(C3)(b)):
 \begin{itemize}
 \item[(C2')] Каждая эксцессная внутренняя вершина $u$ содержится в $\Excess$, и $\tilde e(u)=\{\emptyset\}$.
\end{itemize}

\begin{lemma} \label{lem:after_push}
Предпоток $f$, полученный при достройке, является стабильным и полностью блокирующим.
  \end{lemma}
 \begin{proof}
Рассмотрим ненасыщенный ориентированный путь $P=(u_0,e_1,u_1, \ldots, e_k,u_k)$. Предположим, что либо $u_0\in S$, либо $P$ не доминируется в $u_0$ (когда $u_0$ внутренняя). Тогда ребро $e_1$ закрытое (ввиду~(C1) и (C3)). Значит, вершина $u_1$ подвергалась балансированию и перед этим была эксцессной. Применяя~(C2'),(C3)(b) к вершине $u_1$ и предпотоку $f'$ в момент перед тем балансированием, заключаем, что ребро $e_2$ было закрытым в этот момент. Следовательно, $e_2$ закрытое и для $f$. Рассуждая так и далее, заключаем, что ребра $e_3,\ldots, e_k$ -- тоже закрытые. (Случай $u_k\in T$ невозможен, так как в момент балансирования в вершине $u_{k-1}$ эта вершина была эксцессной, и случай ненасыщенности $e_k$ невозможен.) Так как $e_k$ закрытое, оно должно встречаться в $\deltain(u_k)$ после критического ребра $\hat e(u_k)$ либо совпадать с $\hat e(u_k)$. Тогда из~(C3)(c) следует, что все ребра в $\deltain(u_k)$ после $e_k$ свободны от $f$. Следовательно, $P$ доминируется в $u_k$, и выполняется~\refeq{vk}. Это означает, что $f$ стабильный.

Тот факт, что $f$ полностью блокирующий, выводится из~(C1),(C3),(C2') при аналогичных рассуждениях.
 \end{proof}

Используя указанные выше свойства предпотока, полученного при достройке, можно заключить, что при балансировании на следующей итерации возникает предпоток со свойствами (C1)--(C3). Это обосновывает корректность алгоритма.
 \medskip

\noindent\textbf{Сходимость алгоритма.} Из леммы~\ref{lem:after_push} получаем
  \begin{corollary} \label{cor:stab_flow}
Если на очередной итерации для текущего предпотока $f$ эксцессы всех внутренних вершин становятся нулевыми, то $f$ -- стабильный поток.
  \end{corollary}

В этом случае $f$ -- требуемое решение задачи, и работа алгоритма заканчивается. Поскольку при каждом балансировании или достройке число эксцессных внутренних вершин может как уменьшиться, так и увеличиться, то для установления конечности алгоритма нужен дополнительный анализ.

Для текущего $f$ назовем ребро $e$ \emph{средним}, если оно несвободное и ненасыщенное: $0<f(e)<c(e)$. Пусть $\Gamma=(V_\Gamma,E_\Gamma)$ обозначает подграф в $G$, индуцированный средними ребрами. В частности, $\Gamma$ содержит все несвободные критические ребра. Мы также добавим к $\Gamma$ каждое \emph{свободное} активное ребро (напомним, что активное ребро всегда незакрытое и ненасыщенное). Сделаем следующие наблюдения.
 \smallskip

1) Для ребра $e=uv$ момент его насыщения назовем \emph{событием S}, а момент освобождения (перехода от положительного к нулевому $f(e)$) -- \emph{событием F}. Изменения в $e$ имеют ``однопиковый'' характер: вначале $f(e)$ монотонно увеличивается при достройках в $u$, а при первом уменьшении $f(e)$ ребро $e$ становится закрытым и в дальнейшем $f(e)$ может только уменьшаться (при балансированиях в $v$).
 \smallskip

2) Ребро $e$ может добавиться в граф $\Gamma$ не более двух раз: первый раз  -- когда $e$ становится активным, и второй раз -- когда $e$ становится критическим.
  \smallskip

Пусть $\alpha_S$, $\alpha_F$, $\alpha_M$ обозначают, соответственно, число событий $S$, число событий $F$ и число событий $M$, состоящих в изменении графа $\Gamma$. Из наблюдений выше следует, что

  \begin{numitem1} \label{eq:alpha-Om}
каждое из чисел $\alpha_S,\alpha_F,\alpha_M$ оценивается как $O(m)$.
 \end{numitem1}

Таким образом, для анализа сходимости алгоритма надо оценить число подряд идущих итераций, на которых не происходит событий $S,\,F$ и $M$. Для этого заметим следующее.
\smallskip

3) В процессе алгоритма для каждой внутренней вершины $v$ активное ребро в $\deltaout(v)$ может сдвигаться только вправо (при достройках в $v$), а критическое ребро в $\deltain(v)$ -- только влево (при балансированиях в $v$); при этом статус таких ребер меняется: для старого активного ребра происходит событие $S$ или ребро становится закрытым, а для старого критического ребра -- событие $F$. В силу этого, обозначая $E^+=E^+(f)$ и $E^-=E^-(f)$ множества активных и критических ребер в $\Gamma$, соответственно, получаем, что

 \begin{numitem1} \label{eq:const_Gamma}
в результате операции (балансирования или достройки) граф $\Gamma$ сохраняется тогда и только тогда, когда сохраняются  множества $E^+$ и $E^-$.
 \end{numitem1}

Можно также видеть, что эти множества дают разбиение $E_\Gamma$:
  $$
  E^+\cup E^-=E_\Gamma\quad \mbox{и} \quad E^+\cap E^-=\emptyset
  $$
(учитывая (C3) и то, что критическое ребро $e=uv$ делается закрытым и не может оставаться активным в $\deltaout(u)$). Таким образом, при сохранении $\Gamma$ вместе с сохранением фиксированного разбиения $(E^+,E^-)$ на подряд идущих итерациях каждое изменение предпотока $f$ состоит в его уменьшении в некотором критическом ребре или его увеличении в некотором активном ребре. Заметим также, что активное ребро может стать критическим, но не наоборот. Эти обстоятельства вместе со свойствами~\refeq{alpha-Om} и~\refeq{const_Gamma} позволяют установить следующее

\begin{prop} \label{pr:finite}
Предложенный алгоритм нахождения стабильного потока в сети $N=(G,S,T,c)$ конечный и при целочисленных пропускных способностях $c$ находит целочисленный стабильный поток.
  \end{prop}
\begin{proof}
~Надо показать конечность последовательности итераций, на которых сохраняются $E^+$ и $E^-$. Пусть $\eps$ -- минимальный положительный эксцесс среди внутренних вершин в начале этой последовательности. Предположим по индукции, что перед началом очередного изменения предпотока $f$ на итерации этой последовательности множество $\Excess\cup \New$ эксцессных вершин непусто, и ($\ast$): каждая из них имеет эксцесс не менее $\eps$. Если в этот момент производится балансирование в вершине $v$, то, поскольку $E^-$ сохраняется, балансирование сводится только к уменьшению $f$ на величину $\Delta:=\excess_f(v)\ge\eps$ в критическом ребре $\hat e(v)=uv$. Это обнуляет эксцесс в $v$ и приводит к увеличению эксцесса в вершине $u$ на ту же величину  $\Delta$; следовательно, свойство ($\ast$) верно для нового $f$ (в случае $u\in S$ множество $\Excess\cup\New$ просто уменьшается на один элемент $v$).  Если же производится достройка во внутренней вершине $u$ c активным ребром $\tilde e (u)\ne \emptyset$, то ввиду сохранения $E^+$, достройка сводится к увеличению $f$ в этом ребре $\tilde e(u)=uw$ на величину $\Delta:=\excess(u)\ge \eps$. Это обнуляет эксцесс в $u$ и увеличивает эксцесс в вершине $w$ на величину $\Delta$; следовательно, свойство ($\ast$) верно для нового $f$  (в случае $w\in T$ множество $\Excess\cup\New$ уменьшается на один элемент $u$).

Таким образом, каждое изменение $f$ состоит в его уменьшении на $\ge \eps$ в ребре в $E^-$ либо увеличении на $\ge \eps$ в ребре в $E^+$. Значит, число итераций в данной последовательность не превосходит $c(E^-\cup E^+)/\eps$, что дает конечность алгоритма.

При целочисленном $c$ каждая операция изменяет предпоток в ребре на целое число, и следовательно, результирующий стабильный поток целочисленный.
\end{proof}

Несмотря на конечность алгоритма число итераций в нем может быть очень большим, как показывает следующий пример.
 \medskip

\noindent\textbf{Пример.}  Пусть в $G$ имеются вершины $u,v,w$ и ребра $uv$, $vw$ и $uw$, для которых $uv<_u^+ uw$ и $uw<_w^- vw$. Предположим, ребро $vw$ критическое в $\deltain(w)$, ребро $uv$ критическое в $\deltain(v)$, а ребро $uw$ активное в $\deltaout(u)$. Будем считать, что величины $f(uv)$, $f(vw)$ и $c(uw)-f(uw)$ достаточно большие, эксцесс $\Delta$ в $w$ -- положительный и достаточно малый, а эксцессы в $u$ и $v$ равны нулю. Работа алгоритма может проходить так: сначала $f(vw)$ уменьшается на $\Delta$, затем $f(uv)$ уменьшается на $\Delta$, затем $f(uw)$ увеличивается на $\Delta$. И мы возвращаемся к начальной вершине $u$ c прежним эксцессом $\Delta$ (и при $\excess(u)=\excess(v)=0$) и снова движемся по тому же циклу $w\to v\to u\to w$, и так много раз.
\medskip

В следующем разделе мы расскажем как модифицировать алгоритм, чтобы он приводил к более быстрому решению.

  \section{Модифицированный алгоритм} \label{sec:modif_alg}

  В изложенном выше базовом алгоритме число итераций, на которых граф $\Gamma$ не меняется, может быть очень большим (как показывает Пример в разделе \SEC{bas_alg}); для удобства здесь и далее мы включаем разбиение $(E^+,E^-)$ в описание $\Gamma$. В этом разделе мы описываем модифицированную версию, для которой число изменений предпотока $f$ при постоянном $\Gamma=(V_\Gamma; E^+,E^-)$ имеет порядок $O(n)$. Это эквивалентно тому, что $O(n)$ изменений $f$ приводят к событию $S$, $F$ или $M$ (включая тот случай, когда некоторое активное ребро закрывается и становится критическим).

 Путь в $\Gamma$ назовем \emph{правильным}, если все активные ребра в нем прямые, а все критические ребра обратные. Последовательность итераций с постоянным $\Gamma$ назовем \emph{большой итерацией}. Она начинается с  выбора эксцессной вершины $v_0$ в $(G,f)$, и порядок обработки вершин уточняется следующим образом (полагая, что $v_0$ принадлежит $\Gamma$).
  \smallskip

(A) Если при балансировании очередной эксцессной вершины $v$, состоящем в уменьшении $f$ в критическом ребре $uv$ (которое сохраняется в $\Gamma$), выясняется, что вершина $u$ внутренняя и имеет пустое активное ребро (и тогда все ребра в $\deltaout(u)$ насыщенные или закрытые (ср.~(C3)(b)), ввиду чего достройка в $u$ невозможна), то переходим к балансированию в вершине $u$.
 \smallskip

(B) Если при достройке в вершине $u$, состоящей в увеличении $f$ в активном ребре $e=uw$ выясняется, что вершина $w$ внутренняя и имеет непустое активное ребро $wz$, то далее производится достройка в $w$, а если $\tilde e(w)=\{\emptyset\}$, то переходим к балансированию в $w$. При этом балансировании, если критическим ребром становится то же самое ребро $e=uw$, то оно закрывается и перестает быть активным в $u$, в результате чего граф $\Gamma$ изменяется, и большая итерация заканчивается.
 \smallskip

Из этих уточнений следует, что последовательность обрабатываемых вершин и ребер образует правильный путь $P=(v_0,e_1,v_1,\ldots, e_i,v_i,\ldots)$. Для очередной вершины $v_k$ в $P$ возможны следующие особые ситуации:
 \smallskip

(Q1): $v_k$ совпадает с ранее пройденной вершиной $v_i$;
 \smallskip

(Q2): $v_k$ является терминалом.
 \smallskip

Рассмотрим эти ситуации более подробно. В результате изменения $f$ на $e_1,\ldots, e_k$ эксцессы всех вершин в $P$, кроме $v_k$, обнуляются.
  \smallskip

1) В случае (Q1) выделяем правильный простой цикл $C=(v_i,e_{i+1},\ldots, v_k)$. Подобно тому, как это делается для ротаций в алгоритмах для стабильных $b$-матчингах, распределениях, и др., мы изменяем $f$ вдоль $C$ на максимальную возможную величину $\Delta$, равную минимуму величин $f(e)$ для $e\in E_C\cap E^-$ и $c(e')-f(e')$ для $e'\in E_C\cap E^+$, увеличивая $f$ на $\Delta$ в прямых ребрах и уменьшая на $\Delta$ в обратных ребрах цикла $C$. В результате некоторое ребро в $C$ насыщается или освобождается, и происходит событие $M$ (вместе с $S$ или $F$), завершая большую итерацию. При этом эксцессы всех вершин сохраняются, и можно видеть, что предпоток $f$ продолжает быть стабильным и полностью блокирующим.
 \smallskip

2) В случае (Q2) мы запоминаем $P$ (где теперь эксцессы всех внутренних вершин нулевые) и, в предположении сохранения $\Gamma$, продолжаем большую итерацию, выбирая новую эксцессную вершину $v'_0$ в $(G,f)$ и строя новый правильный путь $P'=(v'_0,e'_1,v_1,\ldots)$.

(a) Если путь $P'$ зацикливается, то действуем как в 1) и заканчиваем большую итерацию.

(b) Если $P'$ встречает предыдущий путь $P$ (ведущий в терминал) в вершине $v'_r=v_i$, то объединяем $P'$ с $P$, получая дерево $\Tscr$ (с корнем в некотором источнике $s\in S$ или ``анти-корнем'' в стоке $t\in T$); при этом все внутренние вершины в $\Tscr$, кроме $v'_r$, имеют нулевой эксцесс. Продолжаем большую итерацию с новой эксцессной вершиной $v''_0$ в $(G,f)$.

(c) Если же $P'$ попадает в терминал отличный от терминала в $P$, то запоминаем $P'$ (помимо $P$) и продолжаем с новой эксцессной вершиной.

(d) В общем случае каждый новый построенный путь, если он не зацикливается и не попадает в новый терминал, встречается либо с деревом $\Tscr$ с корнем в $S$, либо с деревом $\Tscr'$ с анти-корнем в $T$, и мы добавляем этот путь к данному дереву.
  \smallskip

В результате, при сохранении $\Gamma$ (в частности, при отсутствии зацикливания) получаем такую ситуацию:

\begin{numitem1} \label{eq:2trees}
в $\Gamma$ построены несколько непересекающихся правильных деревьев $\Tscr_1,\ldots,\Tscr_k$ с корнями в $S$ и деревьев $\Tscr'_1,\ldots,\Tscr'_{\ell}$ с анти-корнями в $T$, и все внутренние вершины в $G$, не лежащие в этих деревьях, имеют нулевые эксцессы.
  \end{numitem1}

При этом число изменений предпотока $f$ равно $|E_{\Tscr_1}|+\cdots+|E_{\Tscr_k}|+|E_{\Tscr'_1}|+\cdots+|E_{\Tscr'_{\ell}}|\simeq O(n)$. Осталось объяснить как избавляться от этих деревьев.

Предположим $\ell\ne0$, и обозначим $Z$ множество эксцессных вершин в $\Tscr'_1$ (это в точности множество точек ветвления в $\Tscr'_1$). Обойдя $\Tscr'_1$, выделим в нем минимальное поддерево $W$ с анти-корнем в $T$, содержащее $Z$, и перенумеруем ребра в $W$ в топологическом порядке $w_1,\ldots,w_p$ (так, что если $w_j$ лежит на пути, соединяющем $w_i$ с анти-корнем в $\Tscr'_1$, то $i<j$). Перебирая ребра в этом порядке, изменяем $f$ в них естественным образом, получая одно из двух: либо (i) эксцессы всех вершин в $\Tscr'_1$ становятся нулевыми, либо (ii) некоторое промежуточное ребро становится насыщенным или свободным, что заканчивает большую итерацию.

Аналогично действуем и с другими деревьями $\Tscr'_i$ и $\Tscr_j$.

Таким образом, общая работа с деревьями осуществляется за время $O(n)$ и заканчивается либо событием $M$ (вместе с $S$ или $F$), либо избавлением от всех эксцессных вершин в $G$, и следовательно, получением стабильного потока $f$. В течении большой итерации каждое ребро в $\Gamma$ рассматривается не более $O(1)$ раз. Поэтому временн\'ая сложность большой итерации -- $O(n)$.
Этот факт вместе с~\refeq{alpha-Om} приводят к следующему результату.

\begin{prop} \label{pr:modif_alg}
Модифицированный алгоритм находит стабильный поток $f$ в сети $N=(G=(V,E),S,T,c)$ за время $O(nm)$ (причем $f$ целочисленный при целочисленном $c$).
 \end{prop}

\noindent\textbf{Замечание 3.} С достаточной уверенностью можно сказать, что данный алгоритм может быть ускорен до формальной временн\'ой оценки $O(m\,\log n)$ путем применения структур данных как в~\cite{ST1,ST2} для операций на графе $\Gamma$ (таких как выделение и перестройку компонент и циклов в $\Gamma$, агрегированное изменение предпотока на циклах и деревьях, и др.), в том же духе, как подобного рода процедуры для подграфа `` слабых ребер'' ускоряются в алгоритме~\cite{DM}. Мы оставляем проработку этих технических деталей за рамками данной статьи, сохраняя достаточную простоту и практическую применимость изложенных методов.

  \section{Обобщения} \label{sec:general}

Задача о стабильном потоке в сети $N=(G=(V,E),S,T,c)$ может быть обобщена путем введения для каждой внутренней вершины $v$ верхней границы $\gamma(v)\in\Rset_+$ разрешенного эксцесса в $v$. (Как и прежде мы предполагаем условие~\refeq{P}.)
  \medskip

\noindent\textbf{Определение.} Предпоток $f: E\to\Rset_+$ в $N$ назовем $\gamma$-\emph{предпотоком}, если он допустимый (т.е. $f\le c$) и удовлетворяет ограничениям
 \begin{equation} \label{eq:ex-gamma}
0\le \excess_f(v)\le \gamma(v) \quad\mbox{для всех $v\in V-(S\cup T)$.}
  \end{equation}

В прикладных интерпретациях можно понимать ограничение в~\refeq{ex-gamma} как разрешение ``агенту'' $v$ не передавать далее весь продукт, доставленный по входящим ребрам-каналам $e\in\deltain(v)$, а откладывать в запас (или отправлять ``на сторону'') часть полученного продукта в размере, не превышающем $\gamma(v)$.

Мы называем $\gamma$-предпоток \emph{стабильным}, если для каждого ненасыщенного ориентированного $u$--$v$ пути $P$ верно по крайней мере одно из двух: (a) вершина $v$ внутренняя, и $P$ удовлетворяет~\refeq{vk} (при $v_k=v$), т.е. $P$ доминируется в $v$; или (b) вершина $u$ внутренняя, и $P$ \emph{сильно доминируется} в $u$, что означает выполнение~\refeq{v0} (при $v_0=u$) и отсутствие эксцесса в $u$: ~$\excess_f(u)=0$. При $\gamma=0$ это превращается в определение стабильного потока.

Обобщение полученных выше результатов на $\gamma$-предпотоки выглядит так.

\begin{prop} \label{pr:gamma_pre}
Для сети $N=(G=(V,E),S,T,c)$ и вектора $\gamma\in \Rset_+^{V-(S\cup T)}$ стабильный $\gamma$-предпоток существует и может быть найден за время $O(nm)$.
 \end{prop}
 \begin{proof}
~Данная задача сводится к задаче о стабильном потоке в расширенной сети $N'$ с графом $G'=(V,E')$, полученным из $G$ добавлением для каждой внутренней вершины $v$ ребра $vt$ пропускной способности $c(vt):=\gamma(v)$, которое ставится в конец порядка на $\deltaout(v)$, а именно, $e<_v^+ vt$ для всех $e\in\deltaout_G(v)$. Здесь $t$ -- выделенный сток в $T$. Пусть $f'$ -- стабильный поток для $N'$, и $f$ -- его ограничение на $G$. Тогда $f$ -- это $\gamma$-предпоток, и его стабильность следует из стабильности $f'$. (Действительно, если $P$ -- ненасыщенный ориентированный $u$--$v$ путь в $G$, который, будучи рассмотрен как путь в $G'$, доминируется в $u$, то выполняется $f'(ut)=0$, и мы получаем $\excess_f(u)=0$.)
 \end{proof}

Мы можем обобщить задачу далее, вводя для каждой внутренней вершины $v$, помимо $\gamma(v)$ как выше, границу $\beta(v)\in\Rset_+$ для величины $-\excess(v)$. Иначе говоря, мы рассматриваем допустимую (по $c$) функцию $f:E\to\Rset_+$, удовлетворяющую
 \begin{equation} \label{eq:beta-gamma}
-\beta(v)\le\excess_f(v)\le\gamma(v) \quad \mbox{для всех $v\in V-(S\cup T)$.}
 \end{equation}
Такое $f$ назовем $(\beta,\gamma)$-\emph{квазипотоком}. Можно понимать ограничение $-\beta(v)\le \excess_f(v)$ как разрешение ``агенту'' $v$ брать из запаса или привлекать со стороны продукт в размере не более $\beta(v)$ для передачи далее вместе с продуктом, доставленный по входящим ребрам-каналам (при этом, как и прежде, ``агенту'' разрешается отложить в запас или отправить на сторону часть, не превышающую $\gamma(v)$). Проще говоря,  для $(\beta,\gamma)$-квазипотока $f$, если величина $\Delta(v):=\excess_f(v)$ положительная, то ``агент'' $v$ запасает $\Delta(v)$ единиц продукта, а если отрицательная -- берет из запаса $-\Delta(v)$ единиц.

Назовем $(\beta,\gamma)$-квазипоток \emph{стабильным}, если для каждого ненасыщенного ориентированного $u$--$v$ пути $P$ верно по крайней мере одно из двух: (a)  вершина $u$ внутренняя, и $P$ сильно доминируется в $u$ (относительно $\gamma$) в том смысле, что выполненяется~\refeq{v0} (при $v_0=u$), и эксцесс в $u$ неположительный (и не менее $-\beta(u)$); или (b) вершина $v$ внутренняя, и $P$ сильно доминируется в $v$ (относительно $\beta$), что означает выполнение~\refeq{vk} (при $v_k=v$) и неотрицательный эксцесс в $v$ (не более $\gamma(v)$). При $\beta=\gamma=0$ получаем определение стабильного потока.

Для данного обобщения справедливо

\begin{prop} \label{pr:beta-gamma}
Для сети $N=(G=(V,E),S,T,c)$ и векторов $\beta,\gamma\in \Rset_+^{V-(S\cup T)}$ стабильный $(\beta,\gamma)$-квазипоток существует и может быть найден за время $O(nm)$.
 \end{prop}
 \begin{proof}
~Рассмотрим задачу о стабильном потоке в расширенной сети $N'$ с графом $G'=(V',E')$, полученным из $G$: (a) расщеплением каждой внутренней вершины $v$ на две копии $v'$ и $v''$, где $v'$ наследует входящие ребра из $\deltain(v)$, а $v''$ -- выходящие ребра из $\deltaout(v)$ (с теми же пропускными способностями и порядками на них);  (b) добавлением ребра $v'v''$ большой пропускной способности; и (c) добавлении ребра $v't$ пропускной способности $c(v't):=\gamma(v)$ и ребра $sv''$ пропускной способности $c(sv''):=\beta(v)$, где $v't$ и $sv''$ полагаются менее приоритетными, чем ребро $v'v''$, и где $s$ и $t$ -- выделенные источник и сток, соответственно. Пусть $f'$ -- стабильный поток для $N'$, и $f$ -- его ``образ'' в $G$. Можно проверить, что $f$ -- это $(\beta,\gamma)$-квазипоток в $N$, и его стабильность следует из стабильности $f'$. Здесь существеннен тот факт, что для каждой внутренней вершины $v\in V$ по крайней мере одно из значений $f'(sv'')$ и $f'(v't)$ должно быть нулевым (что следует из рассмотрения ненасыщенного пути $(v',v'v'',v'')$).
 \end{proof}

  \section{Дополнительные замечания} \label{sec:concl}

В этом заключительном разделе указываются три дополнительных интересных свойства стабильных потоков и предпотоков. Как и выше, мы рассматриваем ориентированную сеть $N=(G=(V,E),S,T,c)$ c линейными порядками на $\deltain(v)$ и $\deltaout(v)$ для внутренних вершин $v\in V-(S\cup T)$ (и при условии~\refeq{P}).
 \medskip

\noindent\textbf{I.} Для предпотока $f$ в сети $N=(G,S,T,c)$ будем считать величиной $f$ суммарный эксцесс в стоках: $\valu(f):=\excess_f(T)$ ($=\sum_{t\in T} \excess_f(t)$). Справедливо следующее свойство:

 \begin{numitem1} \label{eq:val-pr-f}
для стабильного предпотока $f$ в $N$ найдется стабильный поток $f'$ в $N$ такой, что $f(e)\le f'(e)$ для всех ребер, входящих в $T$, и следовательно, выполняется $\valu(f)\le\valu(f')$.
  \end{numitem1}

Действительно, если $f$ не является потоком (т.е. имеет эксцессную внутреннюю вершину), то $f$ перестраивается в стабильный поток $f'$ применением последовательности итераций базового алгоритма. Каждое балансирование не изменяет значений в ребрах, входящих в $T$, а достройка  может только увеличивать эти значения, что дает требуемое свойство.
 \medskip

\noindent\textbf{II.} ~Что касается величин стабильных потоков, то, обобщая классические результаты для стабильных марьяжей, стабильных двудольных b-матчингов и др., Флейнер~\cite[Sec.~4]{flein} установил, что:

\begin{numitem1} \label{eq:2stab-flow}
для любых стабильных потоков $f$ и $g$ в сети $N$ значения $f(e)$ и $g(e)$ совпадают для каждого ребра $e$, инцидентного терминалу.
 \end{numitem1}

\noindent Как следствие, $\valu(f)=\valu(g)$. Свойство~\refeq{2stab-flow} доказывалось в~\cite{flein} (где рассматривались двухполюсные сети и допускались произвольные ребра, инцидентные терминалам) путем сведения к соответствующему свойству для стабильных распределений. Для наших сетей можно дать прямое и наглядное доказательство.

А именно, свяжем с функцией $f-g$ соответствующее разложение по путям и циклам. Более точно, поскольку $f$ и $g$ не имеют эксцессных внутренних вершин, можно образовать семейство $\Cscr$, состоящее из простых путей и циклов $C$ с весами $\Delta(C)$ так, чтобы выполнялось:
 \smallskip

(i) для каждого $C\in\Cscr$ его прямые ребра $e$ удовлетворяют $f(e)>g(e)$, а обратные ребра $e'$ удовлетворяют $f(e')<g(e')$;
 \smallskip

(ii) для каждого ребра $e\in E$ сумма весов $\Delta(C)$ по путям/циклам $C$, содержащим $e$, равна $|f(e)-g(e)|$;
 \smallskip

(iii) каждый путь в $\Cscr$ соединяет два разных терминала, а каждый цикл содержит не более одного терминала.
 \smallskip

Предположим, $f$ и $g$ не совпадают на некотором ребре $e$, инцидентном терминалу, и рассмотрим случай, когда $e$ выходит из источника $s\in S$ (т.е. $e\in\deltaout (s)$). Для определенности положим $f(e)>g(e)$. Тогда имеется $C\in\Cscr$, содержащее $e$ (как прямое ребро); при этом либо (a) $C$ -- путь из $s$ в $t\in T$, либо (b) $C$ -- путь из $s$ в $s'\in S$ (возможно, $s'=s$).

В случае (a) имеется последовательность $s=v_0,v_1,\ldots, v_k=t$ вершин в $C$ (где $k$ нечетное), для которой часть $C$ от $v_{i-1}$ до $v_i$ имеет только прямые ребра при нечетном $i$, и только обратные ребра при четном $i$. Тогда

\begin{itemize}
\item[($\ast$)] (i) при нечетном $i$ имеется ориентированный путь $P_i$  из $v_{i-1}$ в $v_i$, на ребрах которого $f$ превосходит $g$, и следовательно, $P_i$ ненасыщенный для $g$ и не содержит ребер свободных от $f$; и (ii) при четном $i$ имеется ориентированный путь $Q_i$ из $v_i$ в $v_{i-1}$, на ребрах которого $g$ превосходит $f$, и следовательно, $Q_i$ ненасыщенный для $f$ и не содержит ребер свободных от $g$.
\end{itemize}

\noindent(Такие пути фигурируют в конкатенации $P_1\cdot Q_2^{-1}\cdot P_3\cdot  Q_4^{-1}\cdot\ldots \cdot  Q_{k-1}^{-1}\cdot P_k$ для $C$, где $ Q^{-1}$ обозначает путь, обратный пути $Q$.) Обозначим $p_i$ ($p'_i$) первое (соответственно, последнее) ребро в $P_i$, и обозначим $q_j$ ($q'_j$) первое (соответственно, последнее) ребро в $Q_j$. Можно видеть, что

\begin{itemize}
\item[($\ast\ast$)]
$p'_i,q'_{i+1}\in \deltain (v_i)$ при нечетном $i<k$, и $q_i,p_{i+1}\in \deltaout (v_i)$ при четном $i$.
 \end{itemize}

Теперь применим~($\ast$) и~($\ast\ast$), двигаясь шаг за шагом по последовательности путей $P_1,Q_2,P_3\ldots$\,. Так как $P_1$ начинается в источнике $v_0=s$ и является ненасыщенным для $g$, и так как выполняется $g(q'_2)>0$, то из стабильности $g$ следует, что $q'_2<^-_{v_1} p'_1$. Тогда из стабильности $f$, ненасыщенности $Q_2$ для $f$, неравенства $f(p_3)>0$, и свойства, что $Q_2$ не доминируется в $v_1$ (в силу $f(p'_1)>0$), получаем $p_3<^+_{v_2} q_2$. И так далее. Дойдя до пути $Q_{k-1}$, получаем $p_k<^+_{v_{k-1}} q_{k-1}$. Но тогда последний путь $P_k$ (который заканчивается в стоке $v_k=t$ и является ненасыщенным для $g$) не доминируется в начальной вершине $v_{k-1}$, вопреки стабильности $g$.
 \smallskip

В случае (b) рассуждения схожи. Здесь мы имеем дело с конкатенацией путей вида $P_1,Q_2^{-1},\ldots, P_{k-1},Q_k^{-1}$ (при четном $k$). При этом последний путь $Q_k$ начинается в источнике $v_k=s'$, является ненасыщенным для $f$, и не доминируется в концевой вершине $v_{k-1}$ (в силу $f(p'_{k-1})>0$ и $q'_k<^-_{v_{k-1}} p'_{k-1}$), что противоречит стабильности $f$.

По соображениям симметрии, случаи несовпадений $f$ и $g$ в ребре, входящем в $T$, также невозможны. Этим заканчивается доказательство~\refeq{2stab-flow}.
 \medskip

\noindent\textbf{III.} Предыдущая конструкция может быть продолжена. А именно, рассмотрим два стабильных потока $f$ и $g$ в сети $N=(G,S,T,c)$. Пусть $H$ -- подграф в $G$, порожденный ребрами в $A:=\{e\colon f(e)>g(e)\}$ и $B:=\{e\colon g(e)>f(e)\}$, и назначим ребрам $e$ в $H$ веса $\omega(e):=|f(e)-g(e)|$. Согласно~\refeq{2stab-flow}, $H$ не содержит терминальных вершин. Более того, $\omega$ раскладывается в неотрицательную линейную комбинацию характеристических функций простых правильных циклов (и тем самым $\omega$ может рассматриваться как ``циркуляция'').

Здесь мы говорим, что (не обязательно простой) цикл $C$ в $H$ является \emph{правильным}, если все его прямые ребра принадлежат $A$, а обратные принадлежат $B$. Для вершины $v$ в таком цикле обозначим $\tau(v)$ упорядоченную пару ребер, инцидентных $v$ и следующих в направлении цикла, и назовем вершину $v$  \emph{особой}, если $\tau(v)$ содержит ребро как из $A$, так и из $B$ (эквивалентно, оба ребра в $\tau(v)$ либо входят в $v$, либо выходят из $v$). В этом случае скажем, что пара $\tau(v)=(e,e')$ \emph{левая}, если $e$ предпочтительнее для $v$, чем $e'$; иначе пара называется \emph{правой}. Рассуждая как при доказательстве~\refeq{2stab-flow}, из стабильности $f$ и $g$ можно заключить, что

\begin{numitem1} \label{eq:uniform}
для любого правильного цикла $C$ пары $\tau(v)$  имеют одинаковую ориентацию для всех особых вершин $v$ в $C$: либо все левые, либо все правые.
 \end{numitem1}

Это свойство распространяется на компоненты связности $K$ в $H$:  для всех особых вершин $v$, встречающихся (в правильных циклах) в $K$, пары $\tau(v)$  одновременно либо левые, либо правые (это следует из того, что две такие пары $\tau(v),\tau(v')$ можно включить в один правильный цикл в $K$).
В соответствии с этим выделим четыре типа компонент $K$. А именно, скажем, что $K$ имеет: \emph{тип $A$} (\emph{тип $B$}), если все ребра в $K$ принадлежат $A$ (соответственно, $B$); и имеет \emph{тип $L$} (\emph{тип $R$}), если $K$ имеет ребра в обоих $A$ и $B$ (в этом случае назовем $K$ \emph{богатым}), и ориентации для особых вершин в $K$ левые (соответственно, правые). Эквивалентно,

\begin{numitem1} \label{eq:left}
богатая компонента $K$ имеет тип $L$, если для любой особой вершины $v$ в $K$ в множестве $\deltain(v)$ ребра из $A$ предшествуют ребрам из $B$, а в множестве $\deltaout(v)$ ребра из $B$ предшествуют ребрам из $A$; в случае типа $R$ ситуация противоположная.
 \end{numitem1}

Пользуясь этим, можно явно описать множество стабильных потоков в $N$ в виде решетки (что делается в~\cite{flein} путем сведения к задаче о стабильном распределении и апелляции к соответствующему результату в~\cite{BB}). А именно, определим следующие функции $h,\ell$ на $E$:

\begin{numitem1} \label{eq:h-ell}
$h$ совпадает с $f$ на компонентах типа $A$ и $R$ и совпадает с $g$ на компонентах типа $B$ и $L$, а $\ell$ определяется противоположным образом; на остальных ребрах $e$ полагается $h(e):=\ell(e):=f(e)=g(e)$.
  \end{numitem1}

Можно видеть, что $h$ и $\ell$ являются потоками. А также для каждой внутренней вершины $v$ поток $h$ доминирует потоки $f,g$ на множестве $\deltaout(v)$, и доминируется этими потоками на множестве $\deltain(v)$, в то время как $\ell$ ведет себя противоположным образом. (Здесь для числовых функций $a,b$ на упорядоченном множестве $(S,\prec)$ считаем, что $a$ доминирует $b$, если либо $a=b$, либо имеется $e\in S$ такое, что $a(e)>b(e)$, $a(e')\ge b(e')$ при $e'\prec e$, и $a(e'')\le b(e'')$ при $e\prec e''$.)

Из вышесказанного можно получить (мы опускаем подробности), что $h$ совпадает с указанным в~\cite[Sec.~4]{flein} потоком $f\vee g$, а $\ell$ совпадает с  потокам $f\wedge g$ (в частности, оба $h,\ell$ стабильные).
\bigskip

\noindent\textbf{Благодарность.} Автор выражает благодарность рецензенту за анализ начальной версии статьи и указание работ~\cite{CM,CMS}, неизвестных автору при ее написании.

\end{document}